\documentclass[12pt]{amsart}

\begin{document}
\pagestyle{myheadings}
\markboth{Feigin, Jimbo and Miwa}
{Vertex operator algebra and monomial basis}

\title[Vertex algebras and monomial basis]
{Vertex operator algebra arising from
the minimal series $M(3,p)$ and monomial basis}
\author[B.~L.~Feigin]{Boris~L.~Feigin}
\address{L.~D.~Landau Institute for Theoretical Physics\\
Chernogolovka 142432\\
Russian Federation
}
\email{feigin@feigin.mccme.ru} 
\author[M.~Jimbo]{Michio~Jimbo}
\address{Graduate School of Mathematical Sciences\\
The University of Tokyo, Tokyo 153-8914, Japan}
\email{jimbomic@ms.u-tokyo.ac.jp}
\author[T.~Miwa]{Tetsuji Miwa}
\address{Division of Mathematics, Graduate School of Science\\
Kyoto University, Kyoto 606-8502, Japan}
\email{tetsuji@kusm.kyoto-u.ac.jp}

\dedicatory{Dedicated to Professor B. M. McCoy on the occasion of
his sixtieth birthday}
\date{\today}

\setcounter{footnote}{0}\renewcommand{\thefootnote}{\arabic{footnote}}

\begin{abstract}
We study a vertex operator algebra (VOA) $V$
related to the $M(3,p)$ Virasoro minimal series. 
This VOA reduces in the simplest case $p=4$ to
the level two integrable vacuum module of $\widehat{sl}_2$.
On $V$ there is an action of a commutative current $a(z)$, which is
an analog of the current $e(z)$ of $\widehat{sl}_2$.
Our main concern is the subspace $W$
generated by this action from the highest weight vector of $V$.
Using the Fourier components of $a(z)$, 
we present a monomial basis of $W$ and a semi-infinite monomial basis of 
$V$.
We also give a Gordon type formula for their characters.
\end{abstract}

\maketitle

\vskip20mm

\newpage

\renewcommand{\theequation}{\thesection.\arabic{equation}}
\def\theenumi{\roman{enumi}}
\def\labelenumi{(\theenumi)}
\font\teneufm=eufm10
\font\seveneufm=eufm7
\font\fiveeufm=eufm5
\newfam\eufmfam
\textfont\eufmfam=\teneufm
\scriptfont\eufmfam=\seveneufm
\scriptscriptfont\eufmfam=\fiveeufm
\def\frak#1{{\fam\eufmfam\relax#1}}
\let\goth\frak

\newcommand{\sfrac}[2]{{\textstyle \frac{#1}{#2}}}

\newcommand{\g}{{\goth{g}}}
\newcommand{\slth}{\widehat{\mbox{\twelveeufm sl}}_2} 
\newcommand{\slN}{\mbox{\twelveeufm sl}_N} 
\font\twelveeufm=eufm10 scaled\magstep1
\font\fourteeneufm=eufm10 scaled\magstep2    
\font\seventeeneufm=eufm10 scaled\magstep3   
\font\twentyoneeufm=eufm10 scaled\magstep4  
\newcommand{\slNBig}{\mbox{\seventeeneufm sl}_N} 
\newcommand{\Z}{{\mathbb Z}} 
\newcommand{\C}{{\mathbb C}} 
\newcommand{\R}{{\mathbb R}} 

\newcommand{\F}{{\mathcal F}}
\newcommand{\Ft}{\widetilde{\mathcal F}}
\renewcommand{\H}{{\mathcal H}}
\newcommand{\Cc}{{\mathcal C}}
\newcommand{\Oc}{{\mathcal O}}
\newcommand{\ve}{\varepsilon}
\newcommand{\la}{\lambda}
\newcommand{\Wb}{\overline{W}}
\newcommand{\dz}{\underline{dz}}
\newcommand{\dw}{\underline{dw}}
\newcommand{\dbr}[1]{{\langle\!\langle #1 \rangle\!\rangle}} 
bracket
\newcommand{\G}[3]{{{G}_{x^{#1}}(#2,#3)}}
\newcommand{\dsp}{\displaystyle}
\newcommand{\qi}[1]{{ [  #1 ]_x}}

\newcommand{\nn}{\nonumber}
\newcommand{\bea}{\begin{eqnarray}}
\newcommand{\ena}{\end{eqnarray}}
\newcommand{\be}{\begin{eqnarray*}}
\newcommand{\en}{\end{eqnarray*}}
\newcommand{\lb}[1]{\label{#1}}

\newcommand{\res}{{\mathop{\rm res}}}
\newcommand{\id}{{\rm id}}
\newcommand{\tr}{{\rm tr}}
\newcommand{\bra}[1]{\langle #1 |}        
\newcommand{\ket}[1]{{| #1 \rangle}}      
\newcommand{\br}[1]{{\langle #1 \rangle}}  
\newcommand{\rs}{{r^*}}
\newcommand{\BW}[5]{\Bigl({#1\atop#3}\ {#2\atop#4}\Bigl|#5\Bigr)}

\newcommand{\Remark}{\medskip \noindent {\it Remark.}\quad}
\newcommand{\example}{\medskip \noindent {\it Example.}\quad}

\newtheorem{thm}{Theorem}[section]
\newtheorem{cor}[thm]{Corollary}
\newtheorem{prop}[thm]{Proposition}
\newtheorem{lem}[thm]{Lemma}
\newtheorem{dfn}[thm]{Definition}
\newcommand{\ignore}[1]{}

\setcounter{footnote}{0}
\renewcommand{\thefootnote}{\arabic{footnote})}
\renewcommand{\arraystretch}{1.2}

\setcounter{section}{0}
\setcounter{equation}{0}

\section{Introduction} \lb{sec:1}

Let us recall some basic facts about integrable
$\slth$--modules.
Fix in $\goth{sl}_2$ the standard basis $\{ e, h, f\}$, 
and let $e(z)$, $h(z)$, $f(z)$ be the corresponding currents
of $\slth$:
\be
e(z) = \sum_{i\in\Z} e_i z^{-i-1}, \quad
h(z) = \sum_{i\in\Z} h_i z^{-i-1}, \quad
f(z) = \sum_{i\in\Z} f_i z^{-i-1}.
\en
It is well known that the operators $e(z)^{k+1}$,
$f(z)^{k+1}$ act as zero on any
level $k$ integrable representation of $\slth$.
Conversely,
the relation $e(z)^{k+1}=0$ allows us to describe
integrable representations as follows.

Consider the polynomial ring $A$ in generators
$\{a_i\}_{i \in \Z}$, 
and let $a(z) = \sum_{i\in\Z} a_i z^{-i-1}$.
$A$ is a graded algebra with the assignment $\deg a_i =-i$.
We say that an $A$-module $W$ belongs to the category $\Oc$ 
if it is $\Z$-graded~:~$W = \oplus_{j \in \Z} W_j$,
and $W_j=0$ for sufficiently small $j$. 
On a module $W$ in the category $\Oc$,
certain infinite combinations of the generators $a_i$
have a well-defined action.
For example, each coefficient of
$a(z)^{m}$ ($m=1,2,\cdots$) has a meaning as an operator on $W$.
Denote by $\tilde{A}$ the algebra obtained by
adjoining to $A$ the coefficients of arbitrary polynomials in 
$a(z)$ and its derivatives $\partial a(z)$, 
$\partial^2 a(z),\cdots$ ($\partial=\partial/\partial z$).   
For each non-negative integer $k$,
let $J_{k+1}$ be the ideal generated by the coefficients
of $a(z)^{k+1}$.
We define the quotient algebra $A_{k+1} = \tilde{A}/J_{k+1}$.
It is also convenient to add to $A_{k+1}$ the elements
$D$ and $D^{-1}$ with the relations $D a_iD^{-1} = a_{i+2}$,
$D D^{-1} = D^{-1} D = 1$.

In \cite{FS93}, it was shown that
the algebra $A_{k+1}$ (extended by $D$)
has a remarkable class of irreducible representations
$\pi_{\alpha, \beta}$, where $\alpha$, $\beta$ are
non--negative integers such that $\alpha+\beta = k$.
The space $\pi_{\alpha, \beta}$ has a basis consisting of
semi--infinite monomials $\prod_{i\in\Z} a_i^{b_i}$,
where the exponents $\{ b_i\}_{i\in\Z}$
run over sequences of
non--negative integers satisfying the following conditions.
\begin{description}
\item[(a)] $b_i + b_{i+1} \le k$,
\item[(b)] there exists an $n_0$ such that
$b_n=0$ for $n\le n_0$,
\item[(c)] there exists an $n_1$ such that
for $n\ge n_1$
\be
b_n = \begin{cases}
\alpha & \text{($n$  odd)},\\
\beta  & \text{($n$ even)}.\\
\end{cases}
\en
\end{description}
Actually,  under the identification $e(z) = a(z)$, 
$\pi_{\alpha,\beta}$ is nothing but the irreducible integrable
representation of $\slth$ with highest weight $(\alpha, \beta)$.

At this point, it is natural to ask the following questions.
What will happen if we
replace the relation $a(z)^{k+1}=0$ by something else?
For which kind of relations is it possible to find
an algebra similar to $A_{k+1}$ and
its irreducible representations
which have bases formed by semi--infinite monomials?
These questions are rather complicated,
and 
for `generic' relations we cannot hope for
such a construction.
So it is important to study concrete examples.

Simplest examples arise when we try to generalize
the case of $\slth$ with $k=1$.
Namely, fix some number $m$ and
consider the current $a(z)$ with relations
\be
a(z)^2 = 0, \bigl(\partial a(z)\bigr)^2 =0,
\cdots, \bigl(\partial^m a(z)\bigr)^2 =0.
\en
Such an $a(z)$ can be realized as $a(z)=\phi_{\sqrt{2m}}(z)$,  
where
$\phi_\beta(z)=:e^{\beta\varphi(z)}:$
stands for the vertex operator with momentum $\beta$.
The current $a(z)$ together with
$\partial\varphi(z)$
and
$a^*(z)=\phi_{-\sqrt{2m}}(z)$
form a vertex operator algebra (VOA).
They are analogous to the currents
$e(z)$, $h(z)$, $f(z)$ of $\slth$. 
The $a(z)$ (resp. $a^*(z)$) are commutative,
$\partial\varphi(z)$ generates the Heisenberg algebra,
and $[\partial\varphi(z), a(w)]$
(resp. $[\partial\varphi(z), a^*(w)]$)
have the same form as in the case of $\slth$.
A major difference is in the bracket $[a(z),a^*(w)]$
which becomes a differential polynomial in $\partial\varphi(z)$.

This VOA has $2m$ irreducible representations
enumerated by an integer $l$ with $0\le l <2m$.
These representations 
have a basis consisting of monomials $\prod_{i\in\Z} a_i^{b_i}$,
where
\begin{description}
\item[(a')] $b_i+ b_{i+1} + \dots + b_{i+2m-1} \le 1$,
\item[(b')] There exists an $n_0$ such that
$b_n=0$ for $n\le n_0$,
\item[(c')] There exists an $n_1$ such that
for $n\ge n_1$
\be
b_n = \begin{cases}
1 & \text{($n\equiv l~~\bmod 2m$)},\\
0 & \text{(otherwise)}.\\
\end{cases}
\en
\end{description}

Now let us try to generalize the case of $\slth$ with $k=2$.
The simplest idea is to
replace the relation $a(z)^3=0$ by
two relations $a(z)^3=0$, $a(z)\bigl(\partial a(z)\bigr)^2 =0$.
It is possible to construct a space
which admits an action of $a(z)$ and has
a basis formed by semi-imfinite monomials
$\prod_{i\in\Z} a_{i}^{b_i}$.
The exponents $\{b_i\}$ satisfy conditions similar to
(a'), (b') and (c') above,
wherein the most interesting property (a') is
now replaced by the condition $b_i + b_{i+1} +b_{i+2} \le 2$.

The corresponding VOA can be constructed explicitly.
To do that, recall first the following well--known construction
of $\slth$ with $k=2$.
Consider the VOA obtained as the tensor product of
the Virasoro minimal theory $(3,4)$ (Ising model)
and the lattice vertex algebra generated by $\phi_{\pm1}(z)$.
Let $\psi(z)$ be the $(2,1)$ primary field of
the $(3,4)$ theory (Ising fermion).
Then the formulas
\be
&&a(z) = \psi(z) \phi_1(z),\\
&&a^*(z) = \psi(z) \phi_{-1}(z),
\en
give a realization of the currents $e(z),f(z)$
of $\slth$ at level 2.
In particular $a(z)^3 =a^*(z)^3 =0$.
In the setting above,
let us now
replace the $(3,4)$ theory by the $(3,5)$ theory.
We also replace the Ising fermion
by the $(2,1)$ primary field of the latter,
and the vertex operators $\phi_{\pm1}(z)$
by $\phi_{\pm \sqrt{{3}/{2}}}(z)$.
It turns out that
the resulting current $a(z)$ satisfies the
desired relations
$a(z)^3 =0$ and $a(z)\bigl(\partial a(z)\bigr)^2 =0$.
 
As a next step we can try a larger set of cubic relations,
for example $a(z)^3=0$, $a(z) \bigl(\partial a(z)\bigr)^2 =0$
and $a(z) \partial a(z) \partial^2 a(z) =0$.
But even this case is rather hard to study.
In this paper we try to see what will happen
if we replace the $(3,4)$ or the $(3,5)$ theory
in the previous construction
by some other Virasoro minimal model.
Though the situation in general is obscure,
the above construction goes through for the $(3,p)$ theory.
For example, for the $(3,7)$ theory, 
we get an abelian current $a(z)$ satisfying $5$ cubic relations.

Our main results are the following.
\begin{enumerate}
\item For the $(3,p)$ theory we construct
the current $a(z)$ and describe the set of cubic relations.
\item We find a semi--infinite monomial basis
of the resulting VOA.
\item We obtain a Gordon--type (fermionic) formula for the
characters of VOA and its natural subspace which we call
`principal subspace'.
\end{enumerate}

The text is organized as follows.   
In section 2, after preparing the notation,
we introduce the VOA and the principal subspace,
and state the main results.
Using the `functional model' \cite{FS93}, 
the study of the characters of these spaces is reduced 
to that of certain spaces of symmetric polynomials which arise
as correlation functions
(matrix elements of products of currents). 
In section 3 we determine the structure of three point functions.
We use these results in section 4 to give an upper bound
for the characters.
By comparing it with known fermionic formulas
for the Virasoro minimal characters \cite{BMS},
we find that this bound is in fact exact.
Section 5 is devoted to a combinatorics
to obtain the monomial basis.

\setcounter{section}{1}
\setcounter{equation}{0}
\section{Construction of Vertex Operator Algebras}\lb{sec:2}

\subsection{Notation}\lb{subsec:2.1}
The subject of the present article is a vertex operator algebra
given as a tensor product of two conformal field theories ---
the $(3,p)$ Virasoro minimal series and a free bosonic theory.
First we review a few basic facts about these and fix the notation. 

Let $p\ge 4$ be an integer not divisible by $3$.
The $(3,p)$ minimal series representations
of the Virasoro algebra is characterized by the central charge
\be
c'=1-\frac{2(p-3)^2}{p}.
\en
Let
$T'(z)=\sum_{n\in\Z}L'_n z^{-n-2}$
denote the current of the Virasoro algebra.
Let further
$M_{r,s}=M_{r,s}(3,p)$ ($r=1,2$, $1\le s\le p-1$) be
the irreducible Virasoro module with
central charge $c'$ and highest weight
\bea
\Delta_{r,s}=\frac{(pr-3s)^2-(p-3)^2}{12p},
\lb{delta}
\ena
and let
$\ket{r,s}\in M_{r,s}$ be
the corresponding highest weight vector.
We have $M_{2,s}\simeq M_{1,p-s}$.

To each $M_{r,s}$ there corresponds the
$(r,s)$ primary field $\psi_{r,s}(z)$ of conformal dimension
\eqref{delta}.
We shall particularly be concerned with
the cases $(r,s)=(1,1),(2,1)$.
The corresponding primary fields are the identity $\psi_{1,1}(z)=I$,
and the $(2,1)$ primary field $\psi(z)=\psi_{2,1}(z)$.
The latter is characterized by the operator product expansion
\be
T'(z)\psi(w)=\frac{\Delta_{2,1}}{(z-w)^2}\psi(w)
+\frac{1}{z-w}\partial\psi(w)+O(1),
\en
where $\Delta_{2,1}=(p-2)/4$.
Viewed as an operator
$M_{1+i,1}\rightarrow M_{2-i,1}$ ($i=0,1$), 
$\psi(z)$ has a Fourier mode expansion
\be
\psi(z)=\sum_{n\in\Z}\psi_n z^{-n-\frac{p-2}{2}i}.
\en
The modes are indexed so that $\psi_n$ has degree $-n$:
$[L'_0,\psi_n]=-n\psi_n$.
We normalize $\psi(z)$ by
$\psi_0\ket{1,1}=\ket{2,1}$,
$\psi_0\ket{2,1}=\ket{1,1}$.

Next let $h(z)=\sum_{n\in\Z}h_nz^{-n-1}$ denote the current of the
Heisenberg algebra
\be
[h_m,h_n]=m\delta_{m+n,0}.
\en
For a complex number $\gamma$, let
\be
\F_\gamma=\C[h_{-1},h_{-2},\cdots]\ket{\gamma}
\en
be
the Fock space with the highest weight vector $\ket{\gamma}$ satisfying
\be
h_{n}\ket{\gamma}=0~~(n>0),
\qquad h_0\ket{\gamma}=\gamma\ket{\gamma}.
\en
The Virasoro algebra acts on $\F_\gamma$ through the current
\bea
T''_\lambda(z)=\frac{1}{2}:h(z)^2:+\lambda \partial h(z) 
\lb{T''}
\ena
with the central charge $c''_\lambda=1-12\lambda^2$. 
Here $\lambda$ is a complex number, and the normal ordering rule
$:h_mh_n:=h_mh_n$ ($m\le n$) is implied.
Finally let
$\phi_\alpha(z):\F_{\gamma}\rightarrow \F_{\gamma+\alpha}$ 
denote the chiral vertex operator
\be
\phi_\alpha(z)
=\exp\left(-\alpha\sum_{n< 0}\frac{h_n}{n}z^{-n}\right)
e^{\alpha Q}z^{\alpha h_0}
\exp\left(-\alpha\sum_{n> 0}\frac{h_n}{n}z^{-n}\right),
\en
where $e^{\alpha Q}$ is the isomorphism of vector spaces
$\F_{\gamma}\overset{\sim}{\rightarrow} \F_{\gamma+\alpha}$ such that 
$[e^{\alpha Q},h_n]=0$ ($n\neq 0$), 
$e^{\alpha Q}\ket{\gamma}=\ket{\gamma+\alpha}$.

In the next subsection we consider the Virasoro current
\bea
T(z)=T'(z)+T''_\lambda(z)=\sum_{n\in\Z}L_n z^{-n-2}
\lb{T}
\ena
which acts on $M_{r,s}\otimes \F_{\gamma}$
with the central charge $c=c'+c''_\lambda$.


\subsection{Vertex operator algebra}\lb{subsec:2.2}
The $(2,1)$ field $\psi(z)$ obeys a particularly
simple fusion rule, which is
expressed as an operator product expansion
\bea
&&
(z-w)^{(p-2)/2}\psi(z)\psi(w)
\nn\\
&&\quad
=
I+\frac{p-2}{2c'}(z-w)^2T'(w)+O((z-w)^3)
\lb{phiOPE}\\
&&\quad
=(w-z)^{(p-2)/2}\psi(w)\psi(z).
\nn
\ena
This means that
each matrix element of
$(z-w)^{(p-2)/2}\psi(z)\psi(w)$
can be continued analytically
to $\C^{\times}\times\C^{\times}$,
coincides with that of $(w-z)^{(p-2)/2}\psi(w)\psi(z)$,
and has the expansion \eqref{phiOPE} as $z\rightarrow w$.

Let us modify $\psi(z)$ to get rid of
the factor $(z-w)^{(p-2)/2}$.
Consider the currents
\bea
&&
a(z)=\psi(z)\phi_{\beta}(z),
\lb{a}\\
&&
a^*(z)=\psi(z)\phi_{-\beta}(z),
\lb{a^*}
\ena
where
\be
\beta=\sqrt{\frac{p-2}{2}}.
\en
Note that
\be
\phi_{\beta}(z)\phi_{\pm \beta}(w)
=(z-w)^{\pm(p-2)/2}:\phi_{\beta}(z)\phi_{\pm \beta}(w):.
\en
The currents \eqref{a}, \eqref{a^*}
act on the space $V=\bigoplus_{n\in\Z}V_{n}$
with
\bea
V_{n}=\begin{cases}
M_{1,1}\otimes \F_{n\beta} & \text{($n$ even)},\\
M_{2,1}\otimes \F_{n\beta} & \text{($n$ odd)}.\\
\end{cases}  
\lb{Vn}
\ena
Here $\psi(z)$ acts on the first tensor factor
while $\phi_{\pm\beta}(z)$
acts on the second.
They have Fourier mode expansions in integral powers of $z$,
\be
&&
a(z)=\sum_{n\in \Z}a_n z^{-n-1},
\\
&&
a^*(z)=\sum_{n\in \Z}a^*_n z^{-n-p+3},
\en
where $a_n,a^*_n\in {\rm End}(V)$.
In the case $p=4$,
$V$ is isomorphic to the level $2$ integrable vacuum module of
$\widehat{\goth{sl}_2}$ under the identification
$a(z)=e(z)$, $a^*(z)=f(z)$.

For convenience we choose
\be
\lambda=\beta-\beta^{-1}   
\en
in \eqref{T''}, so that $a(z)$ has conformal dimension $1$
with respect to $T(z)$, \eqref{T}.
On each $V_{n}$ we have an action of $T(z)$ and $h(z)$.  
Under this joint action of the Virasoro and
the Heisenberg algebras, $V_{n}$ is
irreducible and is generated by the vector
\be
v_n
=\begin{cases}
\ket{1,1}\otimes \ket{n\beta}& \text{($n$ even)},\\
\ket{2,1}\otimes \ket{n\beta}& \text{($n$ odd)}.\\
\end{cases}
\en
We call $v_n$ extremal vectors.
They satisfy
\bea
&&
L_m v_n=0,\quad h_m v_n=0, \qquad (m>0),
\lb{high1}
\\
&&
a_m v_{2k+i}=0 \qquad(m\ge -(p-2)k),
\lb{high2}
\\
&&
a^*_m v_{2k+i+1}=0 \qquad(m\ge (p-2)k+2),
\nn\\
&&
a_{-1-(p-2)k}v_{2k+i}=v_{2k+i+1},
\quad
a^*_{1+(p-2)k}v_{2k+i+1}=v_{2k+i},
\lb{high3}
\ena
where $k\in\Z$ and $i=0,1$.
In addition, the extremal vector $v_0$ satisfies
\be
L_0v_0=L_{-1}v_0=0.
\en

Recall that two formal series $A(z)$, $B(z)
\in{\rm End}(V)[[z,z^{-1}]]$ are said to be
mutually local if
$(z-w)^NA(z)B(w)=\varepsilon_{AB} (z-w)^NB(w)A(z)$
holds for some integer $N$ and
$\varepsilon_{AB}=\pm 1$. 
The currents
$T(z)$, $h(z)$, $a(z)$, $a^*(z)$
are pairwise mutually local in this sense.
We have the following operator product expansions as $z\rightarrow w$. 
\be
&&
T(z)T(w)=\frac{c/2}{(z-w)^4}
+\frac{2}{(z-w)^2} T(w)+\frac{1}{z-w}\partial T(w)+O(1),
\\
&&
T(z)h(w)=\frac{-2\lambda}{(z-w)^3}+\frac{1}{(z-w)^2}h(w)
+\frac{1}{z-w}\partial h(w)+O(1),
\\
&&h(z)h(w)=\frac{1}{(z-w)^2}+O(1),
\\
&&T(z)a(w)=\frac{1}{(z-w)^2}a(w)+\frac{1}{z-w}\partial a(w)+O(1),
\\
&&T(z)a^*(w)=\frac{p-3}{(z-w)^2}a^*(w)+\frac{1}{z-w}\partial =
a^*(w)+O(1),
\\
&&
h(z)a(w)=\frac{\beta }{z-w}a(w)+O(1),
\\
&&
h(z)a^*(w)=\frac{-\beta }{z-w}a^*(w)+O(1),
\\
&&
a(z)a(w)=O(1),
\\
&&
a^*(z)a^*(w)=O(1),
\\
&&
a(z)a^*(w)=\frac{1}{(z-w)^{p-2}}\left(I+O(z-w)\right).
\en
Note that there is no singularity
in the expansions of $a(z)a(w)$, $a^*(z)a^*(w)$.
In particular we have
\be
a(z)a(w)=
\left(I+\frac{p-2}{2c'}(z-w)^2T'(w)+O((z-w)^3)\right)
:\phi_\beta(z)\phi_\beta(w):.
\en

\begin{prop}\lb{sec2:prop1}
The space $V$ has the structure of a vertex operator algebra.
\end{prop}
\proof
$V$ is spanned by vectors $P v_0$,
with $P$ running over the
monomials of $\{L_n,h_n,a_n,a^*_n\}_{n\in\Z}$.
Clearly
$L_{-2}v_0,h_{-1}v_0,a_{-1}v_0,a^*_{-p+3}v_0$ are linearly independent.
Hence the statement follows from
the generalities on vertex opeartor algebras
(see e.g. \cite{KacVA}, p.110). 
\qed

The operators $L_0$ and $H=h_0/\beta$
give rise to a $\Z\times\Z$-gradation
$V=\bigoplus_{d,n\in\Z}V_{d,n}$, where
\be
V_{d,n}=\{v\in V\mid L_0 v=dv,~~H v=n v\}.
\en
Quite generally, for any bi-graded vector space
$W=\bigoplus_{d,n\in\Z}W_{d,n}$,
we call
\be
{\rm ch}_{q,z}W={\rm tr}_{W}\left(q^{L_0}z^H\right)
=\sum_{d,n}(\dim\,W_{d,n})\,q^dz^n
\en
the character of $W$.
We shall also consider a graded subspace $U$ of some $W_{n}$.
When there is no fear of confusion, we call
\be
{\rm ch}_qU=\sum_{d}(\dim U_{d,n})q^d
\en
the character of $U$ as well.


\subsection{Principal subspace}\lb{subsec:2.3}
>From now on we will focus attention to the current $a(z)$.
For us it is important to note the following commutativity property.
\begin{prop}
For any $n,m\in\Z$ we have
\bea
&&
[a_n,a_m]=0.
\lb{aa}
\ena
\end{prop}
\proof
This is an immedeate consequence of the
definition \eqref{a} and the
operator product expansion \eqref{phiOPE}.
\qed

Thus the linear span $\goth{n}=\bigoplus_{n\in \Z}\C a_n$
is an abelian Lie subalgebra of $\goth{gl}(V)$.
Let us consider the subspace of $V$
generated over $\goth{n}$ by an extremal vector $v_n$,
\be
V^{n}=U(\goth{n})\,v_n.
\en
We have a sequence of embeddings
\bea
\cdots \longrightarrow V^4
\longrightarrow V^2
\longrightarrow V^0
\longrightarrow V^{-2}
\longrightarrow \cdots.
\lb{emb}
\ena
On the other hand,
the operator $\tau=e^{2\beta Q}$ has the property
\bea
&&
\tau a_n \tau^{-1}=a_{n-p+2}.
\lb{Ta}
\ena
In view of \eqref{high2}, \eqref{high3},
we have an isomorphism of vector spaces
\be
\tau~:~ V^n\overset{\sim}{\longrightarrow} V^{n+2}.
\en
Following \cite{FS93} we call
\be
W=V^0=\C[a_{-1},a_{-2},\cdots]v_0
\en
the principal subspace.

\begin{prop}\lb{sec2:prop2}
The space $V$ coincides with the
inductive limit of the sequence \eqref{emb},
$\cup_N V^{-N}=\varinjlim_{N} V^{-N}$.
\end{prop}
\proof
Set $\tilde{V}=\cup_N V^{-N}$.
We show that
$\tilde{V}$ is invariant under the action of $h_n,L_n$ for all $n$.
The assertion $\tilde{V}=V$ will then follow
from the fact that $v_m\in \tilde{V}$ and the irreducibility of  $V_m$.

Let us verify that
$h_n\tilde{V}\subset\tilde{V}$ ($n\ge N$) holds for any $N$. 
The case $N=0$ is clear from \eqref{high1} and
\bea
&&
[h_n,a(z)]=\beta z^n a(z).
\lb{han}
\ena
Suppose the statement is true for some $N$.
Acting on $v_m$ with the relation
\bea
&&
a(z)^2=:\phi_{\beta}(z)^2:,
\lb{aa1}
\ena
we find that
$h_{N-1}v_m\in\tilde{V}$ for all $m$.
Together with \eqref{han}, this
in turn implies that $h_{N-1}\tilde{V}\subset\tilde{V}$.

Similarly we can show that
$L_n\tilde{V}\subset\tilde{V}$ ($n\ge N$)
for any $N$,
by using 
\bea
&&
[L_n,a(z)]=z^n(z\partial+n+1)a(z), 
\lb{Lan}
\\
&&a(z)\partial^2a(z)=
:\phi_\beta(z)\partial^2\phi_\beta(z):
+
\frac{p-2}{c'}:\phi_{\beta}(z)^2:T'(z).
\lb{aa2}
\ena
\qed


\subsection{Main results}\lb{subsec:2.4}
Since $a_n$'s are commutative, the principal subspace
$W$ is presented as a quotient of the
polynomial ring $A=\C[\xi_{-1},\xi_{-2},\cdots]$ in
indeterminates $\{\xi_{-n}\}_{n\ge 1}$, 
\be
W\simeq A/\goth{a},
\en
where $\goth{a}$ is an ideal of $A$.
The ring $A$ is bi-graded, 
$A=\bigoplus_{d,n\in\Z}A_{d,n}$, 
by assigning the bi-degree $(d,1)$ to $\xi_{-d}$. 
Note that $A_{d,n}=0$ if $d<n$.
We set $A_{n}=\bigoplus_{d\in\Z} A_{d,n}$. 
The following propositions give
some typical elements of $\goth{a}_3=\goth{a}\cap A_3$.

\begin{prop}
We have
\bea
a(z)^2\partial^\nu a(z)=0\qquad (0\le \nu\le p-3).
\lb{difrel}
\ena
Equivalently, the sum
\be
\sum_{i,j,k>0\atop i+j+k=d}
(k-1)(k-2)\cdots (k-\nu)\xi_{-i}\xi_{-j}\xi_{-k}
\en
belongs to $\goth{a}_3$ 
for any $d$ and $\nu$ ($0\le \nu\le p-3$).
\end{prop}
\proof
This is an immediate consequence of
\bea
a(z)^2a(w)=(z-w)^{p-2}\psi(w):\phi_\beta(z)^2\phi_\beta(w):,
\lb{a2a}
\ena 
which follows from \eqref{aa1}.
\qed

In general, $\goth{a}$ contains more complicated elements.
For example, if $p\ge 7$, the relation given below
does not follow from \eqref{difrel} and their derivatives.

\begin{prop}
We have
\be
&&p(p-4)a(z)^2\partial^{p-1}a(z)
+4(p-3)^2 a(z)\partial a(z)\partial^{p-2} a(z)
\\
&&+2(2(p-3)^2-p)a(z)\partial^2a(z)\partial^{p-3}a(z)=0.
\en
\end{prop}
\proof
Using \eqref{a2a} and \eqref{aa2}, we 
expand $a(z)^2a(w)$, $a(z)\partial a(z) a(w)$ and 
$a(z)\partial^2a(z)a(w)$
as $w\rightarrow z$.
Comparing the coefficients of $(z-w)^{p-3}$, we obtain the assertion.
\qed

Our first result states that
all the relations among $\{a_n\}_{n\in\Z}$
are generated from cubic ones.
\begin{thm}\lb{thm1}
The ideal $\goth{a}$ is generated by 
$\goth{a}_{3}=\goth{a}\cap A_3$.
\end{thm}

Proof of Theorem \ref{thm1} will be given in sections 3 and 4.
In the course we shall also find the character of $\goth{a}_3$:
\begin{prop}
\bea
{\rm ch}_q\goth{a}_{3}=
q^3\frac{(1-q^{p-2})(1-q^{p-1})}{(1-q)(1-q^2)(1-q^3)}.
\lb{chA}
\ena
\end{prop}

Our second result concerns the monomial basis of $W$ and $V$.
By the definition, $W$ is spanned by `monomial vectors'
\bea
a_{-\lambda_1}\cdots a_{-\lambda_n}v_0
\lb{mono}
\ena
where $n\ge 0$ and $\lambda_1\ge\cdots\ge \lambda_n\ge 1$. 
We say that \eqref{mono} is admissible if
\bea
\lambda_{i}-\lambda_{i+2}\ge p-2
\qquad
(i=1,\cdots,n-2).
\lb{adm}
\ena

\begin{thm}\lb{thm2}
Admissible monomial vectors \eqref{mono} constitute a
basis of $W$.
\end{thm}
Theorem \ref{thm2} will be proved in section 5.
This theorem leads to the following `semi-infinite'
description of the space $V$ (cf.\cite{FS93,FM99}).
Note that the vector $v_0$ can be formally rewritten as
\be
&&v_0=a_{p-3}^2v_{-2}=a_{p-3}^2a_{2p-5}^2v_{-4}=\cdots\\
&&\quad =a_{p-3}^2a_{2p-5}^2a_{3p-7}^2\cdots v_{-\infty}.
\en
Hence, in view of Proposition \ref{sec2:prop2},
any vector $v\in V$ can be written as a linear combination of the
expressions
\be
v=a_{-m-1}^{\alpha_{-m-1}}a_{-m}^{\alpha_{-m}}a_{-m+1}^{\alpha_{-m+1}}
\cdots v_{-\infty},
\en
where $\{\alpha_j\}_{j=-m-1}^\infty$ runs over
sequences of nonnegative integers subject to the conditions
\be
&&\alpha_j+\alpha_{j-1}+\cdots+\alpha_{j-p+3}\le 2
\qquad \mbox{ for all }j,
\\
&&\alpha_j=\begin{cases}
2& (j\equiv -1~\bmod p-2)\\
0& ({\rm otherwise})\\
\end{cases}
\quad
\mbox{ for sufficiently large $j$ }.
\en

\setcounter{section}{2}
\setcounter{equation}{0}

\section{Three Point Functions}\lb{sec:3}


\subsection{Functional model}\lb{subsec:3.1}
In order to study the structure of $W$,  
we make use of the `functional model' introduced in \cite{FS93} 
for the restricted dual space $W^*$. 
Let us recall this description.
Let
$\Lambda_{n}=\C[x_1,\cdots,x_n]^{S_n}$
denote the space of symmetric polynomials in $n$ variables,
where $S_n$ stands for the symmetric group on $n$ letters.
Let $\Lambda_{d,n}$ be the subspace consisting of
polynomials of homogeneous degree $d$.
The dual vector space of $A_{d,n}$
can be viewed as the space $\Lambda_{d-n,n}$
through the coupling
\be
&&\langle f(x_1,\cdots,x_n), \xi_{-\lambda_1}\cdots =
\xi_{-\lambda_n}\rangle
\\
&&={\rm Res}_{x_1=\cdots=x_n=0}
f(x_1,\cdots,x_n)x_1^{-\lambda_1}\cdots x_n^{-\lambda_n}dx_1\cdots dx_n.
\en
The dual space of $W_{n}\simeq A_{n}/\goth{a}_{n}$ is
the annihilator of $\goth{a}_{n}$, 
\be
I_{n}=
\{f\in \Lambda_{n}\mid \langle f,\xi\rangle=0~~\forall =
\xi\in\goth{a}_{n}\}.
\en
Since the restriction map $V_n^*\rightarrow W_n^*$ is surjective,
the space $I_{n}$ coincides with the set of symmetric polynomials
\be
f_{\psi}(x_1,\cdots,x_n)=
\psi\left(a(x_1)\cdots a(x_n)v_0\right)
\en
where $\psi$ runs over $V_{n}^*$.
In view of the condition \eqref{high1}
and the commutation relations \eqref{han}, \eqref{Lan}, 
the right action of the Heisenberg and the Virasoro algebras on
$V_{n}^*$ translate as
\bea
&&
f_{\psi\cdot h_m}(x_1,\cdots,x_n)=\beta P_m\cdot =
f_{\psi}(x_1,\cdots,x_n)
\qquad (m>0),
\lb{Ward1}
\\
&&
f_{\psi\cdot L_m}(x_1,\cdots,x_n)=l_m(f_{\psi})(x_1,\cdots,x_n)
\lb{Ward2}\\
&&\qquad\qquad\qquad
+(m+1)P_m\cdot f_{\psi}(x_1,\cdots,x_n)
\qquad (m\ge -1),
\nn
\ena
where
\be
P_m=\sum_{j=1}^n x_j^m,\quad
l_m=\sum_{j=1}^n x_j^{m+1}\frac{\partial}{\partial x_j}.
\en
Hence 
$I_{n}$ is an ideal of the polynomial ring
$\Lambda_{n}$, and is
invariant under the `half'
$\bigoplus_{m\ge -1}\C l_m$
of the Virasoro algebra.

Let $\psi_{n}$ be the composition
$V_n\rightarrow \C v_n \rightarrow \C$,
where the first map is the projection along the subspace $\oplus_{d>{\rm =
deg}~v_n}V_{d,n}$.
Set  $\varphi_n=f_{\psi_{n}}$.

\begin{prop}
The ideal $I_{n}$ is generated by polynomials of the form
$l_{m_1}\cdots l_{m_k}(\varphi_n)$ ($m_1,\cdots,m_k\ge -1$).
\end{prop}
\proof
By the Poincar\'e-Birkhoff-Witt theorem, any $\psi\in V^*$ is a linear
combination of elements of the form $\psi_n L_{m_k}\cdots =
L_{m_1}H_{n_l}\cdots H_{n_1}$.
The assertion follows from \eqref{Ward1},\eqref{Ward2}.
\qed

The polynomial $\varphi_n$ is a product of
matrix elements of $\psi(z)$ and $\phi_\beta(z)$,
\bea
&&
\varphi_n(x_1,\cdots,x_n)
\lb{high}\\
&&
=\bra{1+i,1}\psi(x_1)\cdots \psi(x_n)\ket{1,1}
\bra{n\beta}\phi_{\beta}(x_1)\cdots\phi_{\beta}(x_n)\ket{0},
\nn
\ena
where $\bra{r,s}$ denotes the highest
weight vector of the right module
$M_{r,s}^*$,
and $i=0,1$ according to whether $n$ is even or odd.
The second factor is simply $\prod_{i<j}(x_i-x_j)^{\beta^2}$.


\subsection{3 point functions}\lb{subsec:3.2}
Let us study the ideal $I_{3}$ more closely.
For $n=3$, the first factor of \eqref{high} is the four point function 
of the $(2,1)$ field with one point taken to infinity.
As is well known, it satisfies a hypergeometric differential equation.
We thus obtain the following explicit expression
for the polynomial $\varphi_3$.

\begin{prop}\lb{sec3:prop1}
\bea
&&\varphi_3(x_1,x_2,x_3)
\lb{F3}
\\
&&
=
(x_1-x_2)^{p-2}
F\left(1-\frac{p}{3},2-p,2-\frac{2p}{3};\frac{x_3-x_2}{x_1-x_2}\right),
\nn
\ena
where $F(\alpha,\beta,\gamma;z)$ denotes the Gauss hypergeometric 
function.
\end{prop}
Though not apparent, the right hand side of \eqref{F3} is
symmetric in $(x_1,x_2,x_3)$.
\proof
Since $\psi_3L_{-1}=0$ and $\psi_3L_0=(p+1)\psi_3$,
we have $l_{-1}\varphi_3=0$, $l_0\varphi_3=(p-2)\varphi_3$.
Hence we can write
\bea
\varphi_3(x_1,x_2,x_3)
=(x_1-x_2)^{p-2}g\left(\frac{x_3-x_2}{x_1-x_2}\right)
\lb{homog}
\ena
with some polynomial $g(z)$ satisfying $g(0)=1$.
On the other hand, the null vector condition for $\bra{2,1}\in 
M_{2,1}^*$ reads
$\psi_3(pL'_2-3{L'_1}^2)=0$.
Rewriting this equation by using $L_n'=L_n-L_n''$ and
\eqref{Ward1},\eqref{Ward2}, we obtain
\bea
&&\bigl(2p\,l_2-6\,l_1^2+6(p-2)P_1l_1
\lb{hyp}\\
&&\qquad
-(2p-3)(p-2)P_1^2+3(p-2)P_2\bigr)
\varphi_3=0.
\nn
\ena
This implies that, up to a constant multiple,
$g(z)$ is the unique polynomial solution of
\be
z(1-z)\frac{d^2g}{dz^2}
+\left(2-\frac{2p}{3}+\frac{4}{3}(p-3)z\right)\frac{dg}{dz}
-\frac{1}{3}(p-2)(p-3)g=0.
\en
\qed

\begin{prop}\lb{sec3:prop2}
There exist symmetric polynomials $A_1,B_2,Q_{m-1},R_m$ in 
$(x_1,x_2,x_3)$
 such that
\bea
&&
l_1^2(\varphi_3)=A_1 l_1 (\varphi_3)+B_2\varphi_3,
\lb{l12}
\\
&&
l_m(\varphi_3)=Q_{m-1}l_1(\varphi_3)+R_m\varphi_3,
\qquad (m\ge 2).
\lb{lm}
\ena
\end{prop}
\proof
Eq.\eqref{l12} is a consequence of \eqref{hyp} and \eqref{lm} with 
$m=2$.
Let us verify \eqref{lm}. 
Set $z=(x_3-x_2)/(x_1-x_2)$ and note that
\be
&&l_m(z)=Q_{m-1}l_1(z),
\\
&&
Q_{m-1}={\rm Sym}\frac{x_1^{m+1}x_2}{(x_1-x_2)(x_1-x_3)(x_2-x_3)},
\en
where Sym stands for symmetrization.  
Using \eqref{homog}, we obtain
\be
l_m(\varphi_3)-Q_{m-1}l_1(\varphi_3)=(p-2)R_m\varphi_3,
\en
where
\be
&&
R_m=\frac{x_1^{m+1}-x_2^{m+1}}{x_1-x_2}-(x_1+x_2)Q_{m-1}.
\en
It is easy to see that $R_m$ is symmetric.
\qed

\begin{prop}\lb{sec3:prop3}
The ideal $I_{3}$ is generated by $\varphi_3$ and
$l_1(\varphi_3)$,
\be
&&
I_{3}=\Lambda_{3}\varphi_3+\Lambda_{3}l_1(\varphi_3).
\en
Moreover
\be
&&
\Lambda_{3}\varphi_3\cap\Lambda_{3}l_1(\varphi_3)
=\Lambda_{3}(\varphi_3\times l_1(\varphi_3)).
\en
\end{prop}
\proof The first assertion follows from
Proposition \ref{sec3:prop2}.
To see the second, it suffices to show that
$\varphi_3$ and $l_1(\varphi_3)$ do not have a non-trivial common 
factor.
Let
\be
C_n^\nu(z)=\frac{\Gamma(n+2\nu)}{n!\Gamma(2\nu)}
F(2\nu+n,-n,\nu+\frac{1}{2};\frac{1-z}{2})
\en
denote the Gegenbauer polynomial.
The recursion relations
\be
&&(n+2)C^\nu_{n+2}(z)
=2(n+\nu+1)zC^\nu_{n+1}(z)-(n+2\nu)C^\nu_n(z),
\\
&&(1-z^2)\frac{d}{dz}C^\nu_n(z)=-nzC^\nu_n(z)+
(n+2\nu-1)C^\nu_{n-1}(z),
\en
together with $C^\nu_0(z)=1$
imply that $C^\nu_n(z)$ does not have multiple
zeroes provided $2\nu\not\in\Z$.
Our assertion follows from this and \eqref{F3}.
\qed

Proposition \ref{sec3:prop3} enables
us to determine the character of $I_{3}$ as 
\bea
{\rm ch}_q I_{3}
&=&
\frac{q^{p-2}+q^{p-1}-q^{2p-3}}{(1-q)(1-q^2)(1-q^3)}.
\lb{3char}
\ena
By duality we obtain the character  \eqref{chA}
of $\goth{a}_{3}$.

\subsection{Nested structure}\lb{subsec:3.3}
It turns out that there is a nested structure relating
the ideal $I_{3}$
for $p$ with that for  $p-3$.
Here and after, we use the superscript $[p]$ to indicate
the dependence on $p$ of various quantities, e.g.,
$W={W}^{[p]}$, $I_{3}=I_{3}^{[p]}$, $\varphi_3=\varphi_3^{[p]}$.
In the cases $p=4,5$, we define
$I_{3}^{[p-3]}=\Lambda_{3}$,
$\varphi^{[p-3]}_3=1$.

In what follows we set
\be
D_n=\prod_{1\le i<j\le n}(x_i-x_j)^2.
\en

\begin{prop}\lb{sec3:prop4}
We have an exact sequence
\bea
\phantom{DD}
0~\longrightarrow~ D_3\,I_{3}^{[p-3]}
~{\buildrel \iota\over \longrightarrow}~
I_{3}^{[p]}
~{\buildrel \pi\over \longrightarrow}~
(x-y)^{p-2}\C[x,y]
~\longrightarrow ~0,
\lb{exact}
\ena
where
$\iota$ is the inclusion into $\Lambda_{3}$,
and
$\pi$ is the restriction map sending $g(x_1,x_2,x_3)$ to $g(x,y,y)$.
\end{prop}

\proof
First we show that the image of
$\iota$ is contained in $I_{3}^{[p]}$.
In view of Proposition \ref{sec3:prop3} and
\be
&&
D_3 =
l_1(\varphi_3^{[p-3]})=l_1(D_3\varphi_3^{[p-3]})-D_3^{-1}l_1(D_3)\cdot =
D_3\varphi_3^{[p-3]}, 
\en
it suffices to show that $D_3\varphi_3^{[p-3]}$
belongs to $I_3^{[p]}$.

{}From \eqref{F3} we have
\bea
&&l_1(\varphi_3)=(p-2)(x_1-x_2)^{p-1}
\lb{l1phi}\\
&&\quad
\times
\Bigl(
\frac{x_2+x_3}{x_1-x_2}
F(\alpha,3\alpha-1,2\alpha;z)
+(1-z)F(\alpha,3\alpha,2\alpha;z)\Bigr),
\nn
\ena
where $\alpha=1-p/3$ and $z=(x_3-x_2)/(x_1-x_2)$.
Using an identity
\be
&&
\frac{3(3\alpha+1)}{2(2\alpha+1)}z^2(1-z)
F(\alpha+1,3\alpha+2,2\alpha+2;z)
\\
&&
=(z-2)F(\alpha,3\alpha-1,2\alpha;z)
+2(1-z+z^2)F(\alpha,3\alpha,2\alpha;z),
\en
which can be verified via series expansions, we find
\be
k D_3\,\varphi_3^{[p-3]}
=(P_1^2-3P_2)l_1(\varphi_3^{[p]})
+(p-2)(3P_3-P_1P_2)\varphi_3^{[p]},
\en
where $P_j=x_1^j+x_2^j+x_3^j$ and
$k=-9(p-2)(p-4)/2(2p-9)$. 

It is clear that $\pi\circ\iota=0$.
By \eqref{F3} and \eqref{l1phi} we have
$\pi(\varphi_3^{[p]})=(x-y)^{p-2}$ and
$\pi(l_1(\varphi_3^{[p]}))=(p-2)(x+y)(x-y)^{p-2}$. 
Hence $\pi$ maps $I_{3}^{[p]}$ to $(x-y)^{p-2}\C[x,y]$.
Let $(x-y)^{p-2}J$ be the image.
Then $1,x+y\in J$, and $J$ is invariant under multiplication by
$x+2y$ and $2xy+y^2$.
It follows that $J=\C[x,y]$.

Finally we have from \eqref{3char}
\be
{\rm ch}_qI_{3}^{[p]}-q^6{\rm ch}_qI_{3}^{[p-3]}
&=&q^{p-2}\frac{1}{(1-q)^2}\\
&=&
{\rm ch}_q\bigl((x-y)^{p-2}\C[x,y]\bigr),
\en
which shows the exactness of \eqref{exact} at the middle.
\qed


\setcounter{section}{3}
\setcounter{equation}{0}

\section{Gordon Type Filtration}\lb{sec:4}


\subsection{Gordon type filtration}\lb{subsec:4.1}
The goal of this section is to prove Theorem \ref{thm1}.

Let $\widetilde{\goth{a}}^{[p]}$ be the ideal of $A$
generated by $\goth{a}_{3}^{[p]}$.  
Let $\widetilde{I}^{[p]}\subset \Lambda$
($\Lambda=\bigcup_{n\ge 0}\Lambda_n$)
be the annihilator of $\widetilde{\goth{a}}^{[p]}$.
As before, for $p=4,5$
we define $\widetilde{I}^{[p-3]}$ to be $\Lambda$.
We have
\bea
\widetilde{\goth{a}}^{[p]}\subset\goth{a}^{[p]},
\qquad
\widetilde{I}^{[p]}\supset I^{[p]}.
\lb{inc}
\ena
We shall show that in fact the equality holds.

\begin{prop}\lb{sec4:prop1}
Let $m\ge 1$,
$f\in \widetilde{I}_{n}^{[p]}$.
Then
\be
f(x_1,\cdots,x_{n-2m},y_1,y_1,\cdots,y_m,y_m)
\en
is divisible by
\be
\prod_{1\le i\le n-2m\atop 1\le k\le m}(x_i-y_k)^{p-2}
\prod_{1\le k<l\le m}(y_k-y_l)^{2(p-2)}.
\en
\end{prop}
\proof
When $m=1$, this is a consequence of Proposition \ref{sec3:prop4}
and the definition of $\widetilde{I}_{n}^{[p]}$.
The general case follows from the case $m=1$ by further
specializing arguments.
\qed

Consider the map
\be
\Lambda_{n}\longrightarrow \Lambda_{n-2}\otimes\C[y]
\en
which sends $f(x_1,\cdots,x_n)$ to $f(x_1,\cdots,x_{n-2},y,y)$.
The kernel of this map is $D_n\,\Lambda_{n}$.

\begin{prop}\lb{sec4:prop2}
For $n\ge 3$ we have
\be
\widetilde{I}_{n}^{[p]}\cap D_n\,\Lambda_{n}=
D_n\,\widetilde{I}_{n}^{[p-3]}.
\en
\end{prop}

\proof
In the case $n=3$,
we have shown this relation in
Proposition \ref{sec3:prop2}.
Suppose $n\ge 4$.
Then, for $g\in\Lambda_{n}$,
$D_n\, g$ belongs to
$\widetilde{I}_{n}^{[p]}$ if and only if
the coupling $\br{D_n\, g,\xi\eta}$ vanishes for any
$\xi\in\goth{a}_{3}^{[p]}$ and $\eta\in A_{n-3}$.
The latter means that $D_n\, g\in \goth{a}_{3}^{[p]}$
as a polynomial in $(x_1,x_2,x_3)$. 
Hence, by Proposition \ref{sec3:prop2},
this holds if and only if the product
$hg$ with $h=\prod_{1\le i\le 3}\prod_{4\le j\le n}(x_i-x_j)^2$
belongs to the ideal $I_{3}^{[p-3]}$ as a polynomial in $(x_1,x_2,x_3)$.
Suppose $g$ is homogeneous of degree $d$,
and let
\be
&&g=\sum_{k}g'_k(x_1,x_2,x_3) g''_{d-k}(x_4,\cdots,x_n),
\\
&&
h=\sum_{k}h'_k(x_1,x_2,x_3) h''_{6(n-3)-k}(x_4,\cdots,x_n),
\en
be the decomposition into homogeneous components in $(x_1,x_2,x_3)$
with
the indicated degree.
Using $h'_0\neq 0$, we see by induction on $k$
that the above condition holds if and only if
$g'_k\in \widetilde{I}_{3}^{[p-3]}=I_{3}^{[p-3]}$
for all $k$. 
The last condition is equivalent to $g\in
\widetilde{I}_{n}^{[p-3]}$.
\qed

We are now in a position to
estimate the character of $\widetilde{I}^{[p]}$.
For each $n$, consider the following filtration
\bea
F_{-1}=\{0\}\subset F_0\subset \cdots
\subset F_{m-1}\subset F_{m}\subset\cdots
\subset
\widetilde{I}_{n}^{[p]},
\lb{filt}
\ena
where the subspace $F_m$ ($0\le m\le n/2-1$) is defined by
\be
F_m=\{f\in\widetilde{I}_{n}^{[p]}\mid
f\bigl|_{x_{n-2m-1}=x_{n-2m},\cdots,x_{n-1}=x_n}=0
\}.
\en
Then we have the map
\bea
&&F_m/F_{m-1}\longrightarrow \Lambda_{n-2m}\otimes\Lambda_{m},
\lb{map}
\\
&&f(x_1,\cdots,x_n)\mapsto =
f(x_1,\cdots,x_{n-2m},y_1,y_1,\cdots,y_m,y_m),
\nn
\ena
which is injective by the definition.
Combining Propositions \ref{sec4:prop1} with \ref{sec4:prop2} we find
\begin{prop}\lb{sec4:prop3}
The image of the map \eqref{map} is contained in the space
\bea
&&\prod_{1\le i<j\le n-2m}(x_i-x_j)^2
\prod_{1\le i\le n-2m\atop 1\le k\le m}(x_i-y_k)^{p-2}
\prod_{1\le k<l\le m}(y_k-y_l)^{2(p-2)}
\lb{image}\\
&&
\times \widetilde{I}_{n-2m}^{[p-3]}\otimes\Lambda_{m}.
\nn
\ena
\end{prop}

For two formal series $\chi=\sum_{d,n}c_{d,n}q^dz^n$ and
$\chi'=\sum_{d,n}c'_{d,n}q^dz^n$,
we write $\chi\le\chi'$ if $c_{d,n}\le c'_{d,n}$
holds for all $d,n$.
Proposition \ref{sec4:prop3} entails that
\be
{\rm ch}_{q,z}{I}^{[p]}
&\le&{\rm ch}_{q,z}\widetilde{I}^{[p]}
\\
&\le&
\sum_{n,m}\frac{q^{n(n-1)+(p-2)nm+(p-2)m(m-1)}}
{(q)_m}z^{n+2m}
\sum_{d}(\dim\, \widetilde{I}_{d,n}^{[p-3]})\,q^d,
\en
where $(q)_m=\prod_{j=1}^m(1-q^j)$.
Iterating the above inequality $s$ times where
\bea
s=\mbox{ the integre part of  }\frac{p}{3},
\lb{s}
\ena
we are led to the estimate
\bea
&&{\rm ch}_{q,z}I^{[p]}
\lb{chW0}\\
&&\qquad
\le
\sum_{m_0,m_1,\cdots,m_s\ge 0}
\frac{q^{{}^t
{\bf m}
{B}
{\bf m}
+
{}^t
{\bf A}
{\bf m}
}}
{(q)_{m_0}(q)_{m_1}\cdots(q)_{m_s}}
(qz)^{m_0+2m_1+\cdots+2m_s}.
\nn
\ena
Here we have set
\bea
&&
{}^t{\bf m}=(m_0,m_1,\cdots,m_s),
\lb{BA1}
\\
&&B=
\begin{pmatrix}
s                &\frac{p+s-3}{2}&\frac{p+s-4}{2}&\cdots  &\frac{p-2}{2} =
\\
\frac{p+s-3}{2}  &p+s-3          &p+s-4          &\cdots  &p-2\\
\frac{p+s-4}{2}  &p+s-4          &p+s-4          &\cdots  &p-2\\
\vdots           &\vdots         &\vdots         &\ddots  &\vdots\\
\frac{p-2}{2}    &p-2            &p-2            &\cdots  &p-2\\
\end{pmatrix},
\lb{BA2}
\\
&&
{}^t{\bf A}
=(-s,-p+s+1,\cdots,-p+2).
\lb{BA3}
\ena
Since
$I^{[p]}_{d-n,n}\simeq W_{d,n}^*$,  
 \eqref{chW0} gives
an upper bound for the character
${\rm ch}_{q,z}W^{[p]}=
{\rm ch}_{q,qz}I^{[p]}$.

\subsection{Comparison with known characters}
In this subsection we fix $p$ and $s$ as in \eqref{s}, 
and suppress the dependence on $p$ from the notation.
In \cite{BMS}, a large family of `fermionic' formulas
have been obtained for the characters of the Virasoro
minimal series
$M_{r,s}(p',p)$.
The following formula is included as a special case of their result.
\be
{\rm ch}_{q,z}M_{i+1,1}(3,p)=
\sum_{m_0,m_1,\cdots,m_{s-1}\ge 0\atop m_0\equiv i\bmod 2}
\frac{q^{{}^t\overline{\bf m}\overline{B}^{(p)}\overline{\bf m}
+{}^t\overline{\bf A}^{(p)}\overline{\bf m}}}
{(q)_{m_0}(q)_{m_1}\cdots(q)_{m_{s-1}}},
\en
where $i=0,1$, and
if we set $p=3s+\epsilon+1$ ($\epsilon=0,1$) then
\footnote{In eq.(2.21), \cite{BMS}, for the last diagonal entry of $B$
one should apply the rule of the third line rather than the fourth
in the present case.  
We thank B.~M.~McCoy for communicating this to us.}
\be
&&
{}^t\overline{{\bf m}}=(m_0,m_1,\cdots,m_{s-1}),
\\
&&\overline{B}^{(p)}=
\begin{pmatrix}
\frac{s-\epsilon+1}{4}&\frac{s-1}{2}&\frac{s-2}{2}&\cdots  &\frac{1}{2} =
\\
\frac{s-1}{2}  &s-1          &s-2          &\cdots  &1\\
\frac{s-2}{2}  &s-2          &s-2          &\cdots  &1\\
\vdots         &\vdots       &\vdots       &\ddots  &\vdots\\
\frac{1}{2}    &1            &1            &\cdots  &1\\
\end{pmatrix},
\\
&&
{}^t\overline{\bf A}^{(p)}=
(\frac{s+\epsilon-1}{2},s-1,s-2,\cdots,1).
\en
As an immediate consequence we obtain the following  
expression for the character of the space $V$ in
\eqref{Vn}:
\bea
&&{\rm ch}_{q,z}V=
\lb{chV}
\\
&&
\frac{1}{(q)_\infty}
\sum_{m_0,m_1,\cdots,m_{s}\ge 0}
\frac{q^{{}^t{\bf m}{B}{\bf m}+{}^t{\bf A}{\bf m}}}
{(q)_{m_0}(q)_{m_1}\cdots(q)_{m_{s-1}}}
(qz)^{m_0+2m_1+\cdots+2m_s},
\nn
\ena
where the matrix $B$ and the vector ${\bf A}$ are the same as those 
given in
\eqref{BA2}, \eqref{BA3}.

Let us compare the formula \eqref{chV}
with the estimate \eqref{chW0}.
Recall the isomorphism
$\tau:V^{n}\rightarrow V^{n+2}$ in \eqref{Ta}. 
Since
\be
&&
\tau H=(H-2) \tau,
\\
&&\tau (L_0+(p-2)H)=(L_0-2) \tau,
\en
we have
\be
{\rm ch}_{q,z}V^{-2l}
&=&
{\rm tr}_{V^0}
\left(\tau^l q^{L_0}z^{H} \tau^{-l}\right)
\\
&=&
q^{(p-2)l(l+1)}(qz)^{-2l}
\times
{\rm ch}_{q,q^{-(p-2)l}z}W.
\en
Combining this with \eqref{chW0} we find
\bea
&&{\rm ch}_{q,z}V^{-2l}
\lb{chWt}
\\
&&\le
\sum_{m_0,m_1,\cdots,m_s\ge 0}
\frac{q^{{}^t{\bf m}{B}{\bf m}+{}^t{\bf A}{\bf m}}}
{(q)_{m_0}(q)_{m_1}\cdots(q)_{m_s+l}}
(qz)^{m_0+2m_1+\cdots+2m_s}.
\nn
\ena
As $l\rightarrow \infty$, the right hand side tends to
the character \eqref{chV} obtained from the known formulas.
Since $V=\varinjlim_{N}V^{-N}$, we conclude that
the equality holds in the intermediate steps
\eqref{chW0}, \eqref{chWt}.
In summary, we have shown that

\begin{prop}\lb{sec4:prop4}
The map \eqref{map} is an isomorphism.
We have $\widetilde{I}^{[p]}=I^{[p]}$.
\end{prop}

The proof of Theorem \ref{thm1} is now complete.

\bigskip

\noindent{\it Note added.}\quad
The fermionic character formula for $M(3,p)$ is
an example of a general construction due to
G. E. Andrews, ``Multiple series Rogers-Ramanujan type identities'',
Pacific J. Math. 114 (1984) 267--283.
We thank O. Warnaar for drawing attention to this article.


\setcounter{section}{4}
\setcounter{equation}{0}
\font\teneufm=eufm10
\font\seveneufm=eufm7
\font\fiveeufm=eufm5
\newfam\eufmfam
\textfont\eufmfam=\teneufm
\scriptfont\eufmfam=\seveneufm
\scriptscriptfont\eufmfam=\fiveeufm
\def\frak#1{{\fam\eufmfam\relax#1}}
\let\goth\frak

\newcommand{\Bc}{{\mathcal B}}
\newcommand{\Pc}{{\mathcal P}}

\section{Monomial Basis}\lb{sec:5}


\subsection{Spanning set}\lb{subsec:5.1}
The aim of this section is to prove Theorem \ref{thm2}. 
In the present subsection we show first that 
the monomial vectors \eqref{mono} span the space $W^{[p]}$.

Let $\lambda=(\lambda_1,\cdots,\lambda_n)$
($\lambda_1\ge\cdots\ge\lambda_n\ge 0$)
be a partition of length at most $n$.
We set $|\lambda|=\lambda_1+\cdots+\lambda_n$. 
The set of monomial symmetric functions
($S_n^\lambda$ being the stabilizer of $\lambda$)
\be
m_\lambda=\sum_{\sigma\in S_n/S_n^\lambda}
x_{\sigma(1)}^{\lambda_1}\cdots x_{\sigma(n)}^{\lambda_n}
\in\Lambda_{n}
\en
and the set of monomials
\be
\xi_{\lambda}=\xi_{-\lambda_1-1}\cdots \xi_{-\lambda_n-1}
\in A_{n}
\en
are the dual bases to each other.
We introduce a total ordering among partitions of length at
most $n$ by the lexicographic ordering with respect to
$(|\lambda|,\lambda_1,\cdots,\lambda_{n-1})$.
For a symmetric polynomial $f=\sum_{\lambda}c_\lambda =
m_\lambda\in\Lambda_{n}$,
let $\mu$
be the largest element in $\{\lambda\mid c_\lambda\neq 0\}$.
We call $m_{\mu}$ the leading monomial of $f$.

\begin{prop}\lb{sec5:prop1}
Let $\lambda=(\lambda_1,\lambda_2,\lambda_3)$
be admissible in the sense of \eqref{adm}.
Then there exists an element $f\in I_{3}^{[p]}$ whose leading
monomial is $m_\lambda$. 
\end{prop}
\proof
First consider the case $\lambda=(p-2,j,0)$ ($0\le j \le p-2$).
Let
\bea
f_j(x_1,x_2,x_3)=(l_1-(p-2)P_1)^j\varphi_3(x_1,x_2,x_3).
\lb{fj}
\ena
Using the recursion relation
\be
\bigl((z-1)\frac{d}{dz}+\beta\bigr)F(\alpha,\beta,\gamma;z)
=\frac{(\gamma-\alpha)\beta}{\gamma}F(\alpha,\beta+1,\gamma+1;z),
\en
we find
\be
&&f_j(x_1,x_2,0)=(-1)^j 2^{-j}\prod_{i=0}^{j-1}(p-2-i)
\\
&&
\quad
\times
x_1^{p-2}x_2^j
F\left(1-\frac{p}{3},2-p+j,2-\frac{2p}{3}+j;\frac{x_2}{x_1}\right).
\en
This shows that the leading monomial of $f_j$ is 
$m_{(p-2,j,0)}$.

The leading monomial of
$(x_1+x_2+x_3)^{k}(x_1x_2+x_1x_3+x_2x_3)^l
(x_1x_2x_3)^{m} f_{j}$
is $m_{(p-2+k+l+m,j+l+m,m)}$.
The general case is covered by
an appropriate choice of $k,l,m$.
\qed

\begin{prop}\lb{sec5:prop2}
For $n=3$, admissible monomials constitute a basis of
$W_{3}^{[p]}$.
\end{prop}
\proof
For $d\ge 0$,
let $\{\lambda^{(1)},\cdots,\lambda^{(N_d)}\}$
be the set of admissible partitions satisfying $|\lambda|=d$.
By the proposition above, there exist polynomials
$f^{(1)},\cdots,f^{(N_d)}\in I_{3}^{[p]}$ such that
$f^{(j)}=\sum_{i=1}^{N_d} c_{ij}m_{\lambda^{(i)}}$,
where $(c_{ij})_{1\le i,j\le N_d}$
is a triangular matrix with nonzero diagonal entries.
This implies that
the image of the set $X_d=\{\xi_{\lambda^{(i)}}\}_{1\le i\le N_d}$
under the canonical map $A\rightarrow A/\goth{a}$
gives a linearly independent set. 
On the other hand, we have
\be
\sum_{d\ge 0}N_d q^d
&=&\sum_{0\le a\le b\le c\atop c-a\ge p-2}
q^{a+b+c}
\\
&=&\frac{q^{p-2}+q^{p-1}-q^{2p-3}}{(1-q)(1-q^2)(1-q^3)}
\\
&=&
{\rm ch}_qI_{3}^{[p]}.
\en
Therefore $X_d$ is a basis of $I_{d,3}^{[p]}$.
\qed

\begin{prop}\lb{sec5:prop3}
The set of admissible monomials span $W^{[p]}$.
\end{prop}
\proof
Let
$\xi_\lambda=\xi_{-\lambda_1-1}\cdots \xi_{-\lambda_n-1}$
be an arbitrary monomial.
We use the same letter to denote its image in $A/\goth{a}$.
Suppose there is a successive triple $(i-1,i,i+1)$ of indices
for which $(\lambda_{i-1},\lambda_i,\lambda_{i+1})$ is not admissible.
By Proposition \ref{sec5:prop2},
modulo the ideal $\goth{a}_{3}$ the factor
$\xi_{-\lambda_{i-1}-1}\xi_{-\lambda_{i}-1}\xi_{-\lambda_{i+1}-1}$  
can be written as a linear combination of
monomials $\xi_{-\mu_{i-1}-1}\xi_{-\mu_{i}-1}\xi_{-\mu_{i+1}-1}$ such 
that
\be
&&
\mu_{i-1}-\mu_{i+1}\ge p-2>\lambda_{i-1}-\lambda_{i+1},
\\
&&
\mu_{i-1}+\mu_i+\mu_{i+1}
=
\lambda_{i-1}+\lambda_i+\lambda_{i+1}.
\en
With this replacement, $\xi_\lambda$ becomes a linear combination of
monomials $\xi_{\widetilde{\lambda}}$.
We claim that for each such term we have
$\sum_{i=1}^n i\lambda_i>\sum_{i=1}^n i\tilde{\lambda}_i$. 
Since for a given $d$
there are only a finite number of partitions with $|\lambda|=d$,
this procedure terminates after a finite number of steps.
The result is a linear combination of admissible monomials. 

There are four cases to consider.
\be
&&{\rm 1)}:
\mu_{i-1}>\lambda_{i-1},~~\mu_i<\lambda_i,~~\mu_{i+1}>\lambda_{i+1},
\\
&&{\rm 2)}:
\mu_{i-1}>\lambda_{i-1},~~\mu_i<\lambda_i,~~\mu_{i+1}\le\lambda_{i+1},
\\
&&{\rm 3)}:
\mu_{i-1}>\lambda_{i-1},~~\mu_i\ge 
\lambda_i,~~\mu_{i+1}\le\lambda_{i+1},
\\
&&{\rm 4)}:
\mu_{i-1}\le\lambda_{i-1},~~\mu_i>\lambda_i,~~\mu_{i+1}<\lambda_{i+1}.
\en
Let us consider the case 1).
Let $k$ be the number such that $\lambda_{k-1}\ge \mu_{i-1}>\lambda_k$.
Then $\tilde{\lambda}$ has the form
$\tilde{\lambda}_k=\mu_{i-1}$,
$\tilde{\lambda}_j=\lambda_{j-1}$ ($k<j<i$),
$\tilde{\lambda}_{j}=\mu_{j}$ ($j=i,i+1$) and 
$\tilde{\lambda}_{j}=\lambda_{j}$ for the rest.
Then
\be
\sum_{i=1}^n i\lambda_i-\sum_{i=1}^n i\tilde{\lambda}_i
&=&\sum_{j=k}^{i-2}(\mu_{i-1}-\lambda_j)
+\mu_{i-1}-\mu_{i+1}
-\lambda_{i-1}+\lambda_{i+1}
\\
&>&0.
\en
The other cases can be checked similarly.
\qed


\subsection{Counting admissible monomials}\lb{subsec:5.2}
To conclude the proof of Theorem \ref{thm2}, it remains to
show that the number of admissible monomials of a given degree
$(d,n)$ coincides with $\dim W_{d,n}^{[p]}$.
For convenience
we rewrite $\lambda_{i}+1$ in the previous subsection as
$\lambda_i$.
Given $N\ge 1$, let us consider $\lambda=(\lambda_1,\cdots,\lambda_n)$ 
satisfying
\be
&&N\ge\lambda_1\ge\cdots\ge\lambda_n\ge 1,\\
&&\lambda_{i}-\lambda_{i+2}\ge p-2 \qquad (1\le i\le n-2).
\en
Let $\Cc^{[p]}_{N,n}$ denote the set of such partitions.
For $n=0$, $\Cc^{[p]}_{N,0}$ is a single-element set
consisting of the empty partition.
We set $\Cc^{[p]}_{N,d,n}=\{\lambda\in\Cc^{[p]}_{N,n}\mid 
|\lambda|=d\}$.
For a set $X$, we denote its cardinality by $\sharp X$.
In this subsection we derive a recursion relation relating
$\sharp\,\Cc^{[p]}_{N,d,n}$ to  $\sharp\,\Cc^{[p-3]}_{L,d',l}$.

In the sequel we assume that $p\ge 4$.  
Following a similar construction in \cite{FLPW},
we introduce three kinds of transformations $\Bc_1,\Bc_2,\Bc_3$ for
elements of $\Cc^{[p]}_{N,n}$.
The first two are given as follows.
\be
\Bc_1&:&\Cc^{[p-3]}_{L,l}\rightarrow \Cc^{[p]}_{L+2l-2,l},
\\
&&(\lambda_1,\cdots,\lambda_{l-1},\lambda_l)
~\mapsto~(\lambda_1+2(l-1),\cdots,\lambda_{l-1}+2,\lambda_l),
\\
\Bc_2&:&\Cc^{[p]}_{N,n}\rightarrow \Cc^{[p]}_{N+p-2,n+2},
\\
&&(\nu_1,\cdots,\nu_n)~\mapsto~(\nu_1+p-2,\cdots,\nu_n+p-2,1,1).
\en
The third transformation $\Bc_3$ is defined only on a subset of
$\Cc^{[p]}_{N,n}$.
Given an element $\nu=(\nu_1,\cdots,\nu_n)\in\Cc^{[p]}_{N,n}$,
consider the set $S_\nu$
of indices $i$ ($1\le i\le n$) for which one of the following conditions 
hold.
\begin{center}
\begin{tabular}{ll}
($a_i$):& $\nu_i=\nu_{i+1}$,\\
($b_i$):& $\nu_i=\nu_{i+1}+1$ and  $\nu_{i-1}-\nu_{i+1}\ge p-1$,\\
($c_i$):& $\nu_{i-1}>\nu_i>\nu_{i+1}$ and $\nu_{i-1}-\nu_{i+1}= 
p-2$.\\
\end{tabular}
\end{center}
Here we set $\nu_0=+\infty$ and $\nu_{l+1}=-\infty$.
Then $S_\nu$ is empty if and only if
$\nu$ has the form $\Bc_1(\lambda)$ for some $\lambda\in 
\Cc^{[p-3]}_{L,l}$.
For $i\in S_\nu$, set
\be
&&
{\nu'}^{(i)}_j=
\begin{cases}
\nu_j+\delta_{i+1,j}& \mbox{ in the case ($b_i$)},
\\
\nu_j+\delta_{i,j} & \mbox{ otherwise}.
\\
\end{cases}
\en
Let $i_0=\max S_\nu$.
We define $\Bc_3(\nu)$ by the following rule:
If ${\nu'}^{(i_0)}\in\Cc^{[p]}_{N,n}$, then 
$\Bc_3(\nu)={\nu'}^{(i_0)}$.
If ${\nu'}^{(i_0)}\not\in\Cc^{[p]}_{N,n}$, $i_0-1\in S_\nu$ and
${\nu'}^{(i_0-1)}\in\Cc^{[p]}_{N,n}$, then 
$\Bc_3(\nu)={\nu'}^{(i_0-1)}$.
In all other cases $\Bc_3(\nu)$ is not defined.
We have also the inverse transformation of $\Bc_3$.
Notation being as above, let
\be
&&
{\nu''}_j^{(i)}=
\begin{cases}
\nu_j-\delta_{i+1,j}& \mbox{ in the case ($a_i$)},
\\
\nu_j-\delta_{i,j} & \mbox{ otherwise}.
\\
\end{cases}
\en
If ${\nu''}^{(i_0)}\in\Cc^{[p]}_{N,n}$, then we define $\Bc^*_3(\nu)=
{\nu''}^{(i_0)}$. Otherwise $\Bc^*_3(\nu)$ is not defined.
By case checking one can show

\begin{lem}\lb{sec5:lem0}
If $\Bc_3(\nu)$ is defined, then $\Bc_3^*\Bc_3(\nu)$ is defined and 
equals $\nu$.
If $\Bc_3^*(\nu)$ is defined,
then $\Bc_3\Bc_3^*(\nu)$ is defined and equals $\nu$.
\end{lem}

The following assertions can be easily verified.
\begin{lem}\lb{sec5:lem2}
If $S_\nu\neq \emptyset$ and $n\ge 2$, then there exist unique
$k\ge 0$ and $\overline{\nu}\in\Cc^{[p]}_{N-p+2,n-2}$ such that
$\nu=\Bc_3^k\Bc_2(\overline{\nu})$.
\end{lem}
In fact, $k=\max\{j\mid \Bc^{*j}_3(\nu) \mbox{ is defined}\}$.

\begin{lem}\lb{sec5:lem1}
If $\nu=\Bc_2\Bc_1(\lambda)$ with $\lambda\in \Cc^{[p-3]}_{L,l}$,
then $\Bc_3^k(\nu)$ is defined for $0\le k\le 2L-(p-6)(l-2)+2$.
\end{lem}

\begin{lem}\lb{sec5:lem3}
If $\widetilde{\nu}
=\Bc_3^k\Bc_2(\nu)$ is defined for $\nu\in \Cc^{[p]}_{N,n}$,
then $\Bc_3^{\widetilde{k} }\Bc_2(\widetilde{\nu})$ is defined if and 
only if
$0\le \widetilde{k} \le k$.
\end{lem}

In general, consider the composition
\bea
\nu=(\Bc_3^{\mu_m}\Bc_2)\cdots(\Bc_3^{\mu_1}\Bc_2)\Bc_1(\lambda)
\lb{Btr}
\ena
where $\lambda\in\Cc^{[p-3]}_{L,l}$.
It is defined for
\be
0\le \mu_m\le\cdots\le\mu_1\le 2L-(p-6)(l-2)+2
\en
and
\be
\nu\in \Cc^{[p]}_{L+2(l-1)+(p-2)m,l+2m}.
\en
Let $\Pc_{M,m}$ denote the set of partitions
$\mu=(\mu_1,\cdots,\mu_m)$ satisfying
$M\ge\mu_1\ge\cdots\ge\mu_m\ge 0$.
>From the above considerations we obtain a map
\bea
&&\bigcup_{l,m\ge 0\atop l+2m=n}
\Cc^{[p-3]}_{N-2(l-1)-(p-2)m,l}
\times\Pc_{2(N-1)-(p-2)(n-2),m}
\longrightarrow
\Cc^{[p]}_{N,n},
\lb{cmap}
\\
&&
\qquad
(\lambda,\mu)~~\mapsto ~~\nu,
\nn
\ena
where $\nu$ is given by \eqref{Btr}.

Conversely,
using Lemma \ref{sec5:lem2}--\ref{sec5:lem3},
we can put any $\nu\in\Cc^{[p]}_{N,n}$ into the form
\eqref{Btr} in a unique manner.
Hence we have shown

\begin{prop}\lb{sec5:prop4}
The map \eqref{cmap} is a bijection.
\end{prop}

Let us rewrite this result in terms of the generating function
\bea
\widetilde{\chi}^{[p]}_{N,n}(q)
=\sum_{\nu\in\Cc^{[p]}_{N,n}}q^{|\nu|}.
\lb{cchN}
\ena

\begin{prop}\lb{sec5:prop5}
\bea
&&
\widetilde{\chi}^{[p]}_{N,n}(q)
=\sum_{l,m\ge 0\atop l+2m=n}
q^{l(l-1)+(p-2)lm+(p-2)m(m-1)+2m}
\lb{rec1}\\
&&\qquad\times
\Bigl[{2(N-1)-(p-2)(n-2)+m\atop m}\Bigr]
\widetilde{\chi}^{[p-3]}_{N-2(l-1)-(p-2)m,l}(q),
\nn\\
&&
\widetilde{\chi}^{[p]}_{N,n}(q)
=q^n\Bigl[{N-1+n\atop n}\Bigr]
\qquad \mbox{ for $p=1,2$}.
\lb{rec2}
\ena
\end{prop}
\proof
In the bijection \eqref{cmap}, we have the relation
\be
|\nu|=|\lambda|+|\mu|+l(l-1)+(p-2)lm+(p-2)m(m-1)+2m.
\en
Eq.\eqref{rec1} is a consequence of this.
Since
\be
\sum_{\mu\in\Pc_{M,m}}q^{|\mu|}=\Bigl[{M+m\atop m}\Bigr],
\en
and $\Cc^{[p]}_{N,n}=\Pc_{N,n}\backslash \Pc_{N,n-1}$ for $p=1,2$, 
we have for $p=1,2$
\be
\widetilde{\chi}^{[p]}_{N,n}(q)
&=&
\Big[{N+n\atop n}\Bigr]
-\Big[{N-1+n\atop n-1}\Bigr]
\\
&=&
q^n\Big[{N-1+n\atop n}\Bigr].
\en
\qed


\subsection{Proof of Theorem \ref{thm2}}\lb{subsec:5.3}
Recall the filtration \eqref{filt}.
In Proposition \ref{sec4:prop4} we have shown the isomorphism
\bea
&&F_m/F_{m-1}
\nn\\
&&\qquad
\simeq
\prod_{1\le i<j\le l}(x_i-x_j)^2
\prod_{1\le i\le l\atop 1\le k\le m}(x_i-y_k)^{p-2}
\prod_{1\le k<h\le m}(y_k-y_h)^{2(p-2)}
\lb{iso}\\
&&\qquad
\times I^{[p-3]}_{l}\otimes\Lambda_{m},
\nn
\ena
where $l=n-2m$.
Therefore the character
$\chi^{[p]}_{n}(q)={\rm ch}_qW^{[p]}_{n}$
satisfies the recursion relation
\bea
&&
\chi^{[p]}_{n}(q)
=\sum_{l,m\ge 0\atop l+2m=n}
\frac{q^{l(l-1)+(p-2)lm+(p-2)m(m-1)+2m}}{(q)_m}
\chi^{[p-3]}_{l}(q),
\lb{rec3}\\
&&
\chi^{[p]}_{n}(q)
=\frac{q^n}{(q)_n}
\qquad \mbox{ for $p=1,2$}.
\lb{rec4}
\ena

Let us compare this result with
the limit $N\rightarrow\infty$ of \eqref{cchN}
\be
\widetilde{\chi}^{[p]}_{n}(q)=\lim_{N\rightarrow\infty}
\widetilde{\chi}^{[p]}_{N,n}(q).
\en

\begin{prop}\lb{sec5:prop6}
\be
\chi^{[p]}_{n}(q)=\widetilde{\chi}^{[p]}_{n}(q).
\en
\end{prop}
\proof
Proposition \ref{sec5:prop3} shows that
$\chi^{[p]}_{n}(q)\le \widetilde{\chi}^{[p]}_{n}(q)$.
On the other hand, letting $N\rightarrow\infty$ in 
\eqref{rec1},\eqref{rec2},
we find that $\widetilde{\chi}^{[p]}_{n}(q)$
satisfies the same recursion relations
\eqref{rec3},\eqref{rec4} as $\chi^{[p]}_{n}(q)$.
The assertion is proved.
\qed

This completes the proof of Theorem \ref{thm2}.


\subsection{Finitization}\lb{subsec:5.4}
For $N\ge 1$, consider the space
\be
W^{N,[p]}=W^{[p]}/\C[a_{-N-1},a_{-N-2},\cdots]W^{[p]}.
\en
This is a finite dimensional space.
Let $\Lambda^{N}_{n}$ denote
the subspace of $\Lambda_{n}$ consisting of
symmetric polynomials in $n$ variables
whose degree in each variable is at most $N$.
Then the dual vector space of $W^{N,[p]}_{n}$ is
\be
I^{N-1,[p]}_{n}=I^{[p]}_{n}\cap\Lambda^{N-1}_{n}.
\en

In this section we show a `finitized' version of Proposition 
\ref{sec5:prop2},
which states that admissible monomials with
$\lambda_1\le N$ constitute a basis of $W^{N,[p]}_{n}$
(see Theorem \ref{thm3} below).

First let us give a refinement of Proposition \ref{sec3:prop4}.
\begin{prop}\lb{sec5:prop7}
Let
\be
J^{N-1}=\{g\in\C[x,y]\mid {\rm deg}_xg\le N-p+1,~~
{\rm deg}_yg\le 2N-p\}.
\en
Then the sequence
\bea
&&
\lb{exactN}\\
&&
0~\longrightarrow~ D_3\,I_{3}^{N-5,[p-3]}
~{\buildrel \iota\over \longrightarrow}~
I_{3}^{N-1,[p]}
~{\buildrel \pi\over \longrightarrow}~
(x-y)^{p-2}J^{N-1}
~\longrightarrow ~0
\nn
\ena
is exact. Here the maps $\iota,\pi$ are as given in \eqref{exact}.
\end{prop}
\proof
The only nontrivial part is the surjectivity of $\pi$.
We prove it by induction on $M=N-p+1\ge 0$.

Consider the case $M=0$.
Let $f_j(x_1,x_2,x_3)$
be the polynomial \eqref{fj}.
Since $(l_1-(p-2)P_1)\Lambda^{p-2}_{3}\subset \Lambda^{p-2}_{3}$, 
we have $f_j\in I_{3}^{p-2}$.
Using $\pi(\varphi_3)=(x-y)^{p-2}$ we find, 
for $0\le j\le p-2$,
\be
\pi(f_j)&=&\left(x^2\frac{\partial}{\partial x}+
y^2\frac{\partial}{\partial y}
-(p-2)(x+2y)\right)^j((x-y)^{p-2})
\\
&=&\frac{(p-2)!}{(p-2-j)!}(x-y)^{p-2}(-y)^j.
\en
Since $J^{p-2}={\rm span}\{y^j\mid 0\le j\le p-2\}$,
$\pi$ is surjective.

Suppose the statement is true for $M$.
Then there exist polynomials $f_{k,l}\in I^{N-1,[p]}_{3}$
($0\le k\le M, 0\le l\le 2M+p-2$) 
such that $\pi(f_{k,l})=$
\newline
$(x-y)^{p-2}x^ky^l$.
Multiplying $f_{k,l}$
with elementary symmetric polynomials,
we see that
$(x+2y)x^ky^l$, $(2xy+y^2)x^ky^l$ and $xy^2\cdot x^ky^l$
belong to $\pi(I^{N,[p]}_{3})$ 
for $0\le k\le M$, $0\le l\le 2M+p-2$.
Their linear span contains the basis
$\{x^ky^l\mid 0\le k\le M+1,~~0\le l\le 2M+p\}$ of $J^N$.
\qed

Let
\be
&&
\mathcal{A'}^{N,[p]}_{n}=
\{\xi_{-\lambda_1}\cdots \xi_{-\lambda_n}
\mid \lambda=(\lambda_1,\cdots,\lambda_n)\in
\mathcal{C}^{[p]}_{N,n}
\},
\\
&&
\mathcal{A''}^{N,[p]}_{n}
=
\{\xi_{-\lambda_1}\cdots \xi_{-\lambda_n}
\mid \lambda_1> N,~~\lambda_1\ge\cdots\ge\lambda_n\ge 1,
\lambda_i-\lambda_{i+2}\ge p-2\}.
\en
By Theorem \ref{thm2}, the image of $\mathcal{A}^{[p]}_{n}
=\mathcal{A'}^{N,[p]}_{n}\cup\mathcal{A''}^{N,[p]}_{n}$ 
is a basis of $W_{n}^{[p]}$.
Let $\pi^N_{n}:A_n\rightarrow W^{N,[p]}_{n}$ denote the
canonical projection.
Then $W_{n}^{N,[p]}$ is spanned by
$\pi^N_{n}(\mathcal{A'}^{N,[p]}_{n})$.

\begin{prop}\lb{sec5:prop8}
The set $\pi^N_{3}(\mathcal{A'}^{N,[p]}_{3})$
gives a basis of $W_{3}^{N,[p]}$.
\end{prop}
\proof
Set
\be
\chi_{n}^{N,[p]}(q)={\rm ch}_qW^{N,[p]}_{n}.
\en
Proposition \ref{sec5:prop7}
yields the following recursion relation.
\be
&&\chi^{N,[p]}_{3}(q)=q^{p+1}
\frac{1-q^{N-p+2}}{1-q}\frac{1-q^{2N-p+1}}{1-q}
+q^6\chi^{N-4,[p-3]}_{3}(q),
\\
&&
\chi^{N,[p]}_{3}(q)=
q^{3}\Bigl[{N+2\atop 3}\Bigr]
\qquad \mbox{ for $p=1,2$}.
\en
This is the same recursion relation
\eqref{rec1},\eqref{rec2} satisfied by
$\widetilde{\chi}^{[p]}_{N,3}(q)$.
Therefore we have
$\chi_{3}^{N,[p]}(q)=\widetilde{\chi}^{[p]}_{N,3}(q)$,
and the proposition follows.
\qed

\begin{thm}\lb{thm3} For all $N\ge1$ and $n,d\ge 0$ we have
\be
\dim(W^{N,[p]})_{d,n}=\sharp\,\Cc^{[p]}_{N,d,n}.
\en
\end{thm}
\proof
The image of the set
$\mathcal{A''}^{N,[p]}_{n}$ in the subspace
$(\C[a_{-N-1},a_{-N-2},\cdots]W^{[p]})_n$
is linearly independent.
The theorem is equivalent to the statement that
it is also a spanning set.

Let $\lambda=(\lambda_1,\cdots \lambda_n)$ be a partition with
$\lambda_1>N$, and 
consider the process given in the proof of Proposition \ref{sec5:prop3}  
of reducing the monomial $\xi_{-\lambda_1}\cdots\xi_{-\lambda_n}$.
Each step of the reduction consists in replacing a
successive non-admissible triple
$\xi_{-\lambda_{i-1}}\xi_{-\lambda_i}\xi_{-\lambda_{i+1}}$
by a linear combination of admissible ones.
We claim that
the resulting monomials $\xi_{-\mu_1}\cdots \xi_{-\mu_n}$
all belong to the subspace $(\C[a_{-N-1},a_{-N-2},\cdots]W^{[p]})_n$. 
In fact, it is clear in the case $\lambda_{i-1}\le N$.
In the case $\lambda_{i-1}> N$, it is a consequence of
Proposition \ref{sec5:prop8}. 
Repeating this process a finite number of times
we arrive at a linear combination of admissible monomials in
$(\C[a_{-N-1},a_{-N-2},\cdots]W^{[p]})_n$. 
The proof is over.
\qed

Solving the the recursion relation \eqref{rec1}, \eqref{rec2} for
$\widetilde{\chi}_{N,n}(q)$ we are led to the following formula. 

\begin{thm}
\begin{eqnarray*}
{\rm ch}_{q,z}W^{N,[p]}
=\sum_{m_0,m_1,\cdots,m_s\ge 0}
q^{{}^t{\bf m}B{\bf m}+{}^t{\bf A}{\bf m}}
(qz)^{|{\bf m}|}
\prod_{i=0}^s\left[\begin{matrix}P_i+m_i\\m_i\\ \end{matrix}\right],
\end{eqnarray*}
where
${\bf m}={}^t(m_0,m_1,\cdots,m_s)$, $B$ and ${\bf A}$ are 
as in \eqref{BA1}--\eqref{BA3}, and 
\begin{eqnarray*}
&&P_i=
\begin{cases}
N-1-2(B{\bf m})_0-2A_0 & (i=0), \\
2(N-1)-2(B{\bf m})_i-2A_i & (1\le i\le s), \\
\end{cases}
\\
&&|{\bf m}|=m_0+2\sum_{i=1}^s m_i.
\end{eqnarray*}
\end{thm}

\bigskip

\noindent
{\it Acknowledgments.}\quad
M.J. wishes to thank Yaroslav Pugai for explaining him the
content of the work \cite{FLPW},
and Universit\'e de Montpellier II for invitation and warm hospitality
where the final part of this work was done.
This work is partially supported by
the Grant-in-Aid for Scientific Research (B)
no.12440039, Japan Society for the Promotion of Science.


\end{document}